# Optimal droop control placement in distribution network via an exact OPF relaxation method

H. Sekhavatmanesh, *Member, IEEE,* G. Ferrari-Trecate, Senior *Member, IEEE*, and S. Mastellone *Member, IEEE.*

*Abstract—* In the last decade, the integration of Renewable Energy Sources (RES) in distribution networks has been constantly increasing due to their many technical, economical, and environmental benefits. However, the large-scale penetration of RESs is limited by the grid security constraints, e.g., voltage and current limits. The control of inverter-interfaced RESs can guarantee to comply with those constraints while preserving the RESs performance. However, the installation of additional controllable converter units introduces supplementary investment costs and has therefore to be limited. In this paper, a Mixed-Integer Second-Order Cone (MISOCP) optimization problem is developed to optimally place the Q-V and P-V droop controllers for the RES converters. The objectives of the optimization problem are to minimize investment costs, maintenance costs, and the cost of energy purchase from the grid subject to the system constraints. To provide an accurate and convex formulation of these constraints, we adapt an augmented relaxation method, recently proposed to address the optimal power flow problem in radial distribution networks. Our method is tested on a standard IEEE 34-bus network and the results are compared to those provided by alternative approaches.

## I. INTRODUCTION

The ever-increasing penetration of RESs into distribution networks brings several potential technical, economical, and environmental benefits [1]. The optimal design of the RES-integrated distribution networks and their operation are fundamental in maximizing those potential gain. The problem of layout optimization of Distributed Generators (DGs) is addressed via a number of approaches that are reviewed and classified in [1], [2]. Typically, the decision variables of the layout optimization problem are the size and location of the DGs. For RESs such as photovoltaics (PVs), the locations are fixed to places with high solar irradiation levels. However, the injection of RES power generation into the grid without any control might violate grid security constraints such as voltage and current limits. This can be avoided by curtailing the active power generation and controlling the reactive power of RESs. It will require however to install additional controllers with associated converters on RESs, which introduces supplementary investment costs and has therefore to be limited.

In this paper, we determine the optimal location of controllable converter units, while assuming predetermined size. The objective is to minimize the investment and operation costs while preserving the security constraints. The control of DGs can have a centralized or decentralized structure. In the former category, a unique controller, which has the complete knowledge of the network, computes the DG power setpoints [5]. In decentralized architectures, the active and reactive power outputs are adjusted via local controllers using only local voltage measurements [6]. In this paper, we adopt droop controllers as a standard type of local controllers [7]. Most of the literature on droop-controlled DGs focuses on designing the droop curves and optimizing the droop parameters. The placement and sizing of droop-based DGs is typically studied only for islanded Microgrids (MGs), which are reviewed in [8]. Up to our knowledge, the placement of droop controllers in a grid-connected network has not been addressed so far.

The Optimal Power Flow (OPF) approach is the main tool for many operational and investment planning problems in power systems. It ensures electrical security constraints (e.g., voltage and current limits). However, the OPF is inherently a non-convex and NP-hard optimization problem and, consequently, challenging to solve. To overcome this obstacle, five approaches are proposed in the literature, namely, I) Mixed-Integer Linear Optimization approaches, where the AC power flow formulations are either disregarded [9] or approximated in a linear way [10], II) rule-based approaches based on power flow simulations [11], III) non-linear optimization methodologies [12], IV) meta-heuristic approaches [13], and V) convexification methods [14]. The approximations I, II, III, IV might provide a poor quality solution far from the global optimum and in some cases they might even fail to find an existing feasible solution. Moreover, these approaches are typically computationally heavy. In order to define efficient methods for the integration of security constraints into the optimization problem, convex relaxation methods are proposed for the OPF. However, these relaxations are not exact during high RES production periods. Therefore, the obtained solution does not guarantee the electrical security. This is addressed in [15], by augmenting a general cone relaxed formulation with conservative bounds on the voltage and line current magnitudes. The exactness of the resulting Augmented Relaxed OPF (AROPF) formulation is proved under certain conditions that can be checked a priori. This formulation is modified in [16] to account for constant impedance loads in the network, which resulted in a new relaxation method named Modified AROPF (MAROPF).

The problem addressed in this paper is the optimal placement of the droop controllers with associated converters on some pre-installed RESs in a distribution network. The contribution is two folds: first, we model the standard P-V and Q-V droop characteristics as power constant generation units in parallel with constant-impedance loads. As a second contribution, we apply the recently proposed MAROPF formulation to accurately model the security constraints. We will also show that the general cone relaxation method for OPF might lead to the violation of grid security constraints.

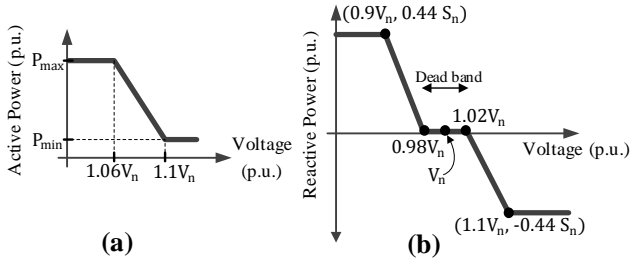

Fig. 1. Standard P-V (a) and Q-V (b) droop characteristics recommended by IEEE 1547.8 [7].

## II. Optimal Placement Problem

An Active Distribution Network (ADN) interconnects consumers and DGs together and to the transmission (or sub-transmission) grid. In this paper, we focus on RESs as DGs and we assume that the ADN is balanced. The active power generated by RESs depends on the environment conditions and cannot be dispatched. This non-dispatchable power injection might violate the security constraints in the grid such as voltage and current limits. Therefore, the penetration level of RESs into ADNs is kept under a specific limit. In this paper, we aim to facilitate high integration of RESs by optimally placing distributed controllers on some of these RESs.

These distributed controllers enable RESs to inject/consume also reactive power as a function of the local voltage measurements. The inverters associated to RESs should have enough capacity for allowing the injection of both active and reactive powers. The size of the inverter for controlling each RES is prespecified according to criteria related to the efficient use of the lifetime of the inverter switches [16]. Therefore, the main decision variable is the placement of controllers and associated converters. The IEEE 1547.8 standard recommends to use a standard characteristic for the adjustment of active and reactive power setpoint, named P-V and Q-V droop [7]. These characteristics together with their standard parameter values are shown in Fig. 1, where $V_n$ is the nominal voltage magnitude. We adopt these droop curves as control rules for the distributed regulators that we aim to optimally place on RESs.

The load and RES power profiles along the whole lifetime of inverters are the input parameters of the optimization problem. These parameters are inherently stochastic, with daily, weekly, seasonally, and yearly variations. A representative number of scenarios should be considered to model these variations. However, integrating all possible scenarios into the optimization problem makes it prohibitive for large-scale systems. An existing solution for dealing with this issue is cluster analysis. According to this method, the original set of scenarios is grouped on the basis of inherent characteristics or similarities [17]. In this paper, we use K-means method to cluster the input parameters into $k$ groups. These parameters include i) the load power profile at each node, ii) the generation power profile of each RES, and iii) the hourly purchasing price of the grid power. The set of cluster centroids constitutes the reduced set of scenarios. Corresponding to each new scenario $s$, there is a parameter $\rho_s$ that represents its probability of occurrence.

### A. Problem Formulation

The notations used for electrical variables are given in Fig. 2 for a typical distribution line. The set of lines is represented with $\mathcal{L}$. Each line is denoted with index $l$, which is also used to refer to the ending node of line $l$. The index $up(l)$ refers to the starting node of line $l$. The slack bus is identified with index 0. The set of nodes hosting RESs is denoted with $\mathcal{G}$. The planning time horizon is $N$ years, which is given as the minimum expected lifetime of the installed inverters. Each year in this period is indexed with $n$. The index $s$ represents a scenario in the reduced set of scenarios $\mathcal{S}$. $\rho_s$ refers to the probability of scenario $s$.

We denote by $v_l$ the square of the voltage magnitude at bus $l$, while the square of the current magnitude in line $l$ is given by $f_l$. The voltage magnitude of the slack bus is fixed to $\sqrt{v_0}$. $v_l^{min}$, $v_l^{max}$, and $I_l^{max}$ denote the squares of the minimum and maximum voltage magnitude limits of bus $l$, and the square of the maximum current magnitude limit of line $l$, respectively. The complex power consumption and power injection at bus $l$ are represented by $s_l$ and $s_l^g$, respectively. The apparent, active, and reactive power flows entering line $l$ are denoted with $P_l$, $Q_l$, and $S_l$, respectively. $z_l$, $r_l$, and $x_l$ represent, respectively, the longitudinal impedance, resistance, and reactance of line $l$.

The investment planning problem is formulated in the following as a Mixed Integer Second Order Cone Programming (MISOCP). The main decision variable is $x_l$, which is defined for each RES at node $l \in \mathcal{G}$ as a binary variable, indicating whether at this node a controller is installed or not. All the other variables are indexed with $s$, which represents their values for scenario $s \in \mathcal{S}$.

$$\underset{x_l}{\text{Minimize}}\ F^{obj} = w_{inv}.F^{inv} + w_{op}.F^{op} \quad (1)$$

$$F^{inv} = \sum_{l \in \mathcal{G}} x_l \left( C_l^F + C_l^{inv} S_l^{max} \right) + \sum_{n=1}^{N} \left( \frac{1}{(1+\alpha)^{n-1}} \sum_{l \in \mathcal{G}} x_l C_l^M \right) \quad (2)$$

$$F^{op} = \sum_{s \in \mathcal{S}} \rho_s C_s^{grid} P_{1,s}^{pur} \quad (3)$$

subject to:

$$\begin{cases} 0 \leq P_{1,s}^{cns} \\ P_{1,s} \leq P_{1,s}^{cns} \end{cases} \quad (4)$$

$$\sum_{l \in \mathcal{G}} x_l \leq X^{max} \quad (5)$$

***Droop control model*** (6)-(8)

***Power flow model*** (9)-(14)

As per (1), the objective function comprises investment ($F^{inv}$) and operational ($F^{op}$) terms. Equation (2) describes the investment costs through three terms, namely, i) the capital cost of installing a controller and the accompanying inverter at node $l$ ($C_l^F$), ii) the investment costs of them at node $l$ ($C_i^{inv}$), which is related to the inverter size $S_l^{max}$, and iii) the maintenance cost of the inverter at node $l$ ($C_l^M$) brought to the year of the investment. The annual interest rate is denoted by $\alpha > 0$. The operational objective term is formulated in (3) as the expected value of the total cost of power purchase from the upper grid. This power is equal to the power flow in the substation line ($P_{1,s}$), when positive.

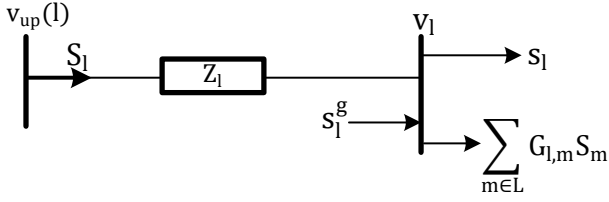

Fig. 2. The model of a distribution line with the notations used for the formulation of the power flow equations

In order to obtain a linear model of the grid power purchase, an auxiliary variable $P_{1,s}^{pur}$ is introduced, which is constrained in (4) [18]. $C_s^{grid}$ is the cost of purchasing 1 p.u. power from the grid in scenario $s$. Constraint (5) limits the total number of controllers to be installed on RESs to a maximum value ($X^{max}$). This represents the limitation on the budget of DNOs for installing the controllers and the required inverters. The other constraints are defined in the following of this section.

*1) Droop Control Model*

The authors have already developed in [16], a linear model for the droop characteristics shown in Fig. 1. In this model, the active and reactive power generation are expressed as a linear function of the square of the measured local voltage. Here we integrate this model into the proposed optimization problem. Constraints (6) and (7) express the active power generation of each RES at node $l$ and scenario $s$, being $M$ a large positive coefficient. In these formulations, $v_l^{0p}$ and $\alpha_l^p$ are functions of the droop parameters indicated in Fig. 1, and $p_{l,s}^{ava}$ is the available active power of the RES at node $l$ and scenario $s$.

$$v_l^{0p} - M(1 - y_{l,s}) \le v_{l,s} \le v_l^{0p} + My_{l,s} \quad (6)$$

$$\begin{cases} -M \cdot y_{l,t} \le p_{l,s}^g - p_{l,s}^{ava} \le 0 & (7.a) \\ -M(2 - y_{l,s} - x_l) \le p_{l,s}^g - \left(p_{l,s}^{ava} - \alpha_l^p(v_{l,s} - v_l^{0p})\right) \\ \le M(2 - y_{l,s} - x_l) & (7.b) \end{cases}$$

Constraint (6) enforces $y_{l,t}$ to be equal to one when $v_{l,s}$ is larger than $v_l^{0p}$ and to zero when $v_{l,s}$ is lower than $v_l^{0p}$. If $y_{l,s} = 0$, (7.a) is imposed and (7.b) is weakened to a constraint which always holds. It means that $p_{l,s}^g$ is forced to be qual to $p_{l,s}^{ava}$. The P-V droop controller should be activated when $y_{l,s} = 1$ and $x_l = 1$. In this case, (7.a) is weakened to an always holding constraints and (7.b) forces $p_{l,s}^g$ to be equal to $\left(p_{l,s}^{ava} - \alpha_l^p(v_{l,s} - v_l^{0p})\right)$. According to [16], this expression provides an approximation of the standard P-V curve shown in Fig. 1.a. For this, $v_l^{0p}$ and $\alpha_l^p$, used in (6) and (7), should be equal to the expressions given in the following:

$$v_i^{0p} = 2\sqrt{v_0 * 1.06V_n} - v_0$$

$$\alpha_i^p = -\frac{P_i^{max} - P_i^{min}}{0.04V_n * 2\sqrt{v_0}}$$

where, $v_0$ is the square of the initial guess of the voltage magnitude, used in the Taylor approximation proposed in [16].

The reactive power generation of each RES for each scenario $s$ is described by the following inequalities.

$$\begin{cases} -M \cdot x_l \le q_{l,s}^g \le M \cdot x_l & (8.a) \\ -M(1 - x_l) \le q_{l,s}^g - \left(-\alpha_l^q(v_{l,s} - v_l^{0q})\right) \le M(1 - x_l) & (8.b) \end{cases}$$

If $x_l = 0$, (8.a) forces $q_{l,s}^g$ to be zero, while (8.b) is relaxed and does not put any constraint on $q_{l,s}^g$. If $x_l = 1$, (8.a) is relaxed and (8.b) forces $q_{l,s}^g$ to be equal to $\left(q_l^{g0} - \alpha_l^q(v_{l,s} - v_l^{0q})\right)$. As proved in [16], This expression provides a Taylor approximation of the standard Q-V droop characteristic shown in Fig. 1, if $\alpha_l^q$ and $v_l^{0q}$ are set as follows:

$$v_i^{0q} = 2\sqrt{v_0 * V_n} - v_0$$

$$\alpha_i^q = -\frac{0.44 S_i^{max}}{0.1 V_n * 2\sqrt{v_0}}$$

*2) Power Flow Model*

Constraints (9)-(12), provided hereafter, define the set of power flow formulations for each scenario $s$, and for each bus (or line) $l$ (if not mentioned otherwise). We apply the MAROPF relaxation method to ensure that the security limits are satisfied if the droop controllers are installed on the selected RESs. The authors have proved in [16] that the MAROPF guarantees the feasibility of the solution in presence of linear P-V and Q-V droop controllers. The following constraints express the power flow formulation according to the general relaxation method [14].

$$\begin{cases} s_{l,t}^g = p_{l,t}^g + jq_{l,t}^g & (9.a) \\ v_{l,t} = v_{up(l),t} - 2\Re(z_l^* S_{l,t}) + |z_l|^2 f_{l,t} & (9.b) \\ S_{l,t} = s_{l,t} - s_{l,t}^g + \sum_{m \in \mathcal{L}} G_{l,m} S_{m,t} + Z_l f_{l,t} & (9.c) \\ f_{l,t} \ge \dfrac{|S_{l,t}|^2}{v_{l,t}} & (9.d) \end{cases}$$

The formulation of the active ($p_{l,s}^g$) and reactive ($q_{l,s}^g$) power generation of RESs were developed in (7) and (8), respectively. The nodal voltage equation is provided in (9.b) for each bus $l$. Constraint (9.c) describes the active-/reactive power balances at the end buses of each line $l$. Constraint (9.d) is the relaxed version of the current flow equation. During overvoltage conditions, where maximum voltage limit is binding, the optimal solution of this relaxation method might not satisfy the original constraint (i.e. the equality condition in (9.c)). In this case, the following additional constraints should be added according to the MAROPF method [16]. Following [16], we impose security limits on a set of auxiliary variables ($\bar{f}, \hat{S}, \underline{S}, \bar{S}, \hat{v}$) to ensure the feasibility of the solution. Variables $\bar{f}$ and $\hat{v}$ represent upper bounds on line current and nodal voltage magnitudes. $\hat{S} = \hat{P} + j\hat{Q}$ and $\hat{v}$ are defined in (11) according to the DistFlow equations, where the line current flows are assumed zero. The formulation of auxiliary variables $\underline{S} = \underline{P} + j\underline{Q}$ are introduced in (12.a) as lower bounds for $S = P + jQ$. Finally, auxiliary variables $\bar{S} = \bar{P} + j\bar{Q}$ and $\bar{f}$ are defined in (12.b) and (12.c) so that they are upper bounds respectively for $S = P + jQ$ and $f$ (see [16] for further details). Together with (9), the following constraints provide the MAROPF formulation, where $(s^g, S, v, f, \hat{S}, \hat{v}, \underline{S}, \bar{S}, \bar{f})$ are the optimization variables.

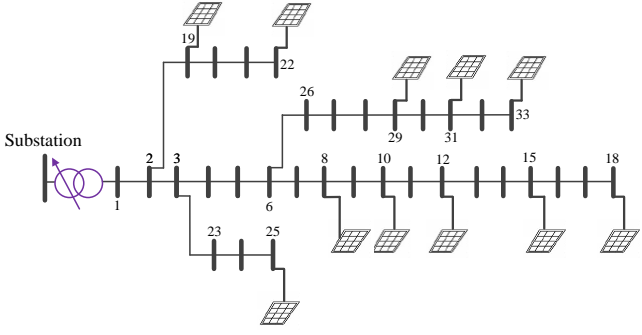

Fig. 3. One line diagram of the modified 34-bus distribution grid [19].

$$\begin{cases} M \cdot y_{i,t} \leq \hat{p}_{l,t}^g - p_{i,t}^{ava} \leq 0 & (10.a) \\ -M(2 - y_{i,t} - x_i) \leq \hat{p}_{l,t}^g - \left(p_{i,t}^{ava} - \alpha_i^p(\hat{v}_{i,t} - v_i^{0p})\right) & (10.b) \\ \hat{p}_{l,t}^g - \left(p_{i,t}^{ava} - \alpha_i^p(\hat{v}_{i,t} - v_i^{0p})\right) \leq M(2 - y_{i,t} - x_i) & (10.c) \\ -M \cdot x_i \leq \hat{q}_{l,t}^g \leq M \cdot x_i & (10.c) \\ -M(1 - x_i) \leq \hat{q}_{l,t}^g - \left(-\alpha_i^q(\hat{v}_{i,t} - v_i^{0q})\right) \leq M(1 - x_i) & (10.c) \end{cases}$$

$$\begin{cases} \hat{s}_{l,t}^g = \hat{p}_{l,t}^g + j\hat{q}_{l,t}^g & (11.a) \\ \hat{S}_{l,t} = s_{l,t} - \hat{s}_{l,t}^g + \sum_{m \in \mathcal{L}} G_{l,m}\hat{S}_{m,t} & (11.b) \\ \hat{v}_{l,t} = \hat{v}_{up(l),t} - 2\Re(z_l^*\hat{S}_{l,t}) & (11.c) \end{cases}$$

$$\begin{cases} \underline{S_{l,t}} = s_{l,t} - s_{l,t}^g + \sum_{m \in \mathcal{L}} G_{l,m}\underline{S_{m,t}} & (12.a) \\ \overline{S_{l,t}} = s_{l,t} - s_{l,t}^g + \sum_{m \in \mathcal{L}} G_{l,m}\overline{S_{m,t}} + Z_l\overline{f_{l,t}} & (12.b) \\ \max\left\{|\overline{P_{l,t}}|^2, |\underline{P_{l,t}}|^2\right\} + \max\left\{|\overline{Q_{l,t}}|^2, |\underline{Q_{l,t}}|^2\right\} \leq \overline{f_{l,t}}v_{up(l)} & (12.c) \end{cases}$$

$$v_l \geq v_l^{min}, \quad \hat{v}_l \leq v_l^{max} \tag{13}$$

$$\bar{f}_l \leq I_l^{max} \tag{14}$$

Assuming that the power loss term $(Z_l f_l)$ is zero, (11.b) and (11.c) express Kirchhoff's current and voltage laws, respectively. The formulation of $\hat{p}_{l,s}^g$ and $\hat{q}_{l,s}^g$ are given in (10), which are similar to (7) and (8), except that instead of $v$, the variable $\hat{v}$ is used to formulate the droop characteristic. As proved in [16], these yield that $\hat{v}$ is an upper bound for $v$. The lower bounds on the active and reactive power flows are defined in (12.a), which is similar to the formulation of (11.b), except that instead of $\hat{s}_{l,t}^g$, the variable $s_{l,t}^g$ is used. For defining upper bounds on $P$ and $Q$, the inequality (12.b) is added, which is similar to the Kirchhoff's current low in (9.c), while instead of $f$, its upper bound $\bar{f}$ is used. The auxiliary variable $\bar{f}$ is defined in (12.c). To account for the reverse current flow in the grid, the maximum of the absolute values of upper- and lower bounds on power flow variables are used in (12.c) to define $\bar{f}$. Note that all the upper and lower bounds are independent from $f$, which is a basic requirement in the MAROPF formulation. According to the defined upper bounds on $v$ and $f$, the voltage and current limits will be accounted for by using the conservative constraints given in (13) and (14), respectively.

### III. SIMULATION RESULTS

The proposed MISOCP formulation for the investment planning of distributed controllers is applied on a 34-bus network, shown in Fig. 3. This 12.66kV distribution network is a standard IEEE test benchmark, which is connected to the sub-transmission grid through a substation. The voltage of this substation is assumed fixed to 1.05 p.u., which is realized in practice using a tap changing transformer, as shown in Fig. 3. The detailed nodal and branch data is given in [19].

The test network is modified by assuming 11 PV resources on different buses. The lifetime of the inverters and distributed controllers is assumed to be 5 years. The investment planning is computed for this period, while considering the variation and uncertainty of the load and PV generation powers. The hourly load profile at different nodes of the network are according to the industrial, commercial, residential, and rural load patterns reported in [20]. The hourly forecast data for the power generation of each PV resource is obtained according to the real data reported in [21]. Total PV capacity is 9.85 MW at full insolation periods in the day, which is 17.63 times larger than the total load powers. As mentioned in section II, the apparent power capacity of the inverters to be placed on PVs is optimized to achieve the most efficient functionality of the inverter switches. This capacity for each PV source is assumed to be 110 % of its maximum active power generation capacity.

For the price of power purchase from the grid, typical average price of electricity in a region in the southern part of Switzerland is used [22]. As mentioned in section II, the variation and uncertainty scenarios of each load profile, each PV generation profile, and grid power price profiles are reduced using K-means method. To cover daily (morning, midday, evening), weekly (weekday and weekends), and seasonally variations of parameters, the whole data set is clustered into 24 (=3*2*4) groups. The centroid of each cluster is regarded as a scenario and integrated into the optimization problem. The under- and overvoltage limits and the overcurrent limit are, respectively, set to 0.90, 1.10, and 7.5 p.u. The scaled weighting factors $w_{inv}$ and $w_{op}$ are set to 0.8 and 1.0 p.u., respectively. It is assumed that maximum 6 number of distributed controllers are allowed to be installed ($X^{max} = 6$). In this simulation study, the developed MISOCP problem is solved using the Branch-and-Bound method. The optimization model is implemented on a PC with an Intel(R) Core i7 CPU and 6 GB RAM; and solved in Matlab/Yalmip environment, using Gurobi solver.

Fig. 4 shows the voltage and current profiles in the network, when there is no distributed controller is installed. These profiles correspond to the maximum voltage and current magnitudes among all the considered scenarios $s \in S$. As it can be seen, the upper voltage and current limits are violated due to the high penetration of PV sources into the network.

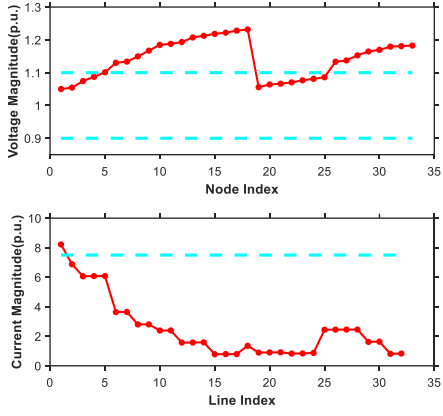

Fig. 4. Maximum voltage and current profiles without any control on the PV resources. Dashed lines represent limits.

Table I. Numerical results of simulation cases 1, 2, and 3.

| case | Optimal controller places | Objective function (p.u.) | Maximum voltage (p.u.) | Maximum current (p.u.) | Computation Time (s) |
|---|---|---|---|---|---|
| 1 | {10, 12, 15, 18, 31, 33} | 0.4375 | 1.095 | 6.620 | 25.18 |
| 2 | {12, 15, 18, 31, 33} | 0.3741 | 1.141 | 8.951 | 20.36 |
| 3 | -- | 0 | 1.232 | 8.220 | 1.772 |

Three simulation case studies are compared in this section. In case 1, the developed formulation is applied on the test system. This is compared with case 2 to highlight the impact of integrating the developed linear model of the P-V and Q-V droop characteristics into the optimization problem. In case 2, the same simulation conditions as those in case 1 are applied, while excluding the model of the droop controllers from the optimization formulation (constraints (7), (8), and (10)). In order to present the benefit of integrating MAROPF formulation into the optimization problem, case 1 is compared with case 3, where instead of the MAROPF relaxation, the general cone relaxation method is used for the power flow formulation [14].

Table I shows the obtained optimal solution of the investment planning in each simulation cases. The last column shows the computation time required for each case. As seen, the computational load of the developed formulation is sufficiently low, which confirms its scalability to real applications. According to the results of case 1, six distributed controllers should be installed. However, as shown in Table I, the number of installed controllers obtained in cases 2 and 3 are smaller and therefore, the obtained objective values are also smaller. To check the voltage and current security limits at each simulation case, the droop controllers are positioned on the model of the test network, according to the obtained optimization solution. This model is implemented in Matlab/Simulink environment. Load flow simulations are carried out for each variation scenario, to obtain the real voltage and current profiles. The fourth and fifth columns of Table I gives the maximum voltage and current magnitudes out of these profiles for each test case. The results show that the overvoltage and overcurrent limits are all respected in case 1, whereas in cases 2 and 3, they are violated at some nodes and lines of the network. It means that the obtained solutions in cases 2 and 3 are not feasible, although they resulted in lower objective values with respect to case 1. The reason for this infeasibility in case 2 is that the droop control characteristics are not modeled in the optimization problem.

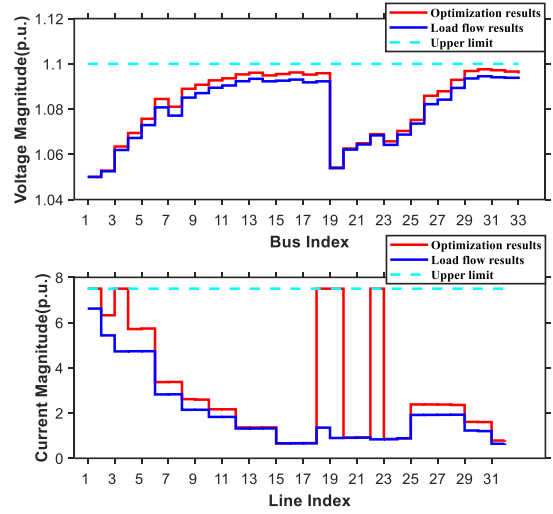

Fig. 5. The maximum voltage and current results in case 1

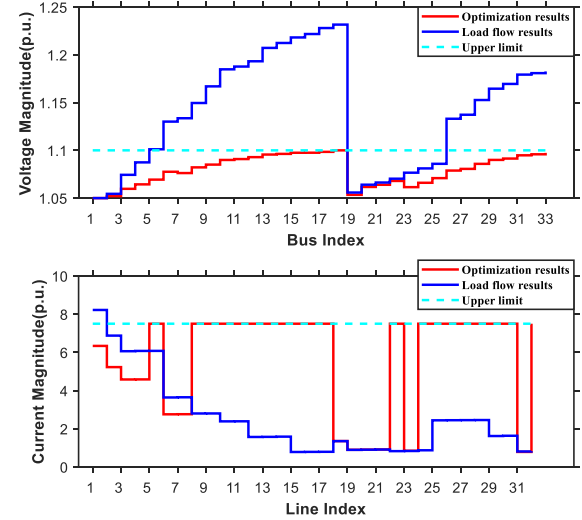

Fig. 6. The maximum voltage and current results in case 3

To better understand the reason for the infeasibility in case 3, the maximum voltage and current profiles obtained in case 1, and 3 are shown in Fig. 5 and Fig. 6, respectively. According to Fig. 5, the real voltage and current values are not matching with the values obtained from the optimization problem. However, they are all below the maximum limits, meaning that the developed optimization problem guarantees the feasibility of the solution. However, Fig. 6 shows that the mismatch in case 3 between the optimization and real voltage and current profiles leads to the violation of the maximum limits. This mismatch is due to the huge reverse power flow in the network that causes the upper voltage limits to bind and, consequently, leading to an inexact solution of the general cone relaxation method.

## IV. Conclusion

In this paper, a MISOCP formulation is developed for optimally placing distributed controllers and accompanying inverters on RESs in a distribution network. As the first contribution, the standard droop control characteristic recommended by IEEE 1547.8 is modelled and integrated into the optimization problem in a linear fashion. The second contribution was to accurately model the electrical security constraints using a convex formulation. In this regard, a recently-published relaxation method for AC power flow equations, named MAROPF, is applied and integrated into the investment planning.

The functionality of the proposed investment planning approach is illustrated on a standard 34-bus test benchmark. Using a posteriori power flow simulation, it is shown that the optimization solution satisfies all the security constraints. The performance of the proposed optimization layout approach is compared to those of standard approaches for the 34-bus test benchmark. In one case, the model of the droop control characteristic is removed from the optimization problem. In the second case, the MAROPF formulation is replaced by the general cone relaxation of AC power flow equations. Both cases resulted in infeasible solutions, violating the upper voltage and upper current limits. Finally, the short computation time reported for the developed approach shows that it can be used for real applications with large-scale networks.